\documentclass[12pt]{iopart}
\usepackage{amssymb}
\usepackage{amsfonts}
\usepackage[english]{babel}
\usepackage{euscript}
\usepackage{bbm}
\usepackage{xfrac}
\usepackage{mathrsfs}
\usepackage{amsbsy}

\begin{document}

\title[The point equivalence problem for ODEs of the second order]{%
The point equivalence problem for ordinary differential
equations of the second order}
\author{Oleg I. Morozov${}^{*}$}

\begin{abstract}
We use \'E. Cartan's method  to solve the problem of equivalence of the second order ordinary differential
equations with respect to the pseudogroup of point trans\-for\-ma\-ti\-ons.
\end{abstract}

\vfill
\noindent
\footnotesize{*~The work was partially supported by the Joint Grant 09-01-92438-KE\_a of RFFR (Russia) and
Consortium E.I.N.S.T.E.IN (Italy)}
%\hfill
%\today

\maketitle

%%%%%%%%%%%%%%%%%%%%%%%%%%%%%%%%%%%%%%%%%%%%%%%%%%%%%%%%%%%%%%%%%%%%%%%%%%%%%%%%%%%%%%%%%%%%%%%%%%%%%%%%%%

\section{Introduction}
The problem of finding the necessay and sufficient conditions of equivalence of ordinary differential equations
of the second order
\begin{equation}
u_{xx}= F(x,u,u_x)
\label{main}
\end{equation}
with respect to the pseudogroup of local diffeomorphisms
\begin{equation}
\widetilde{x}= \varphi(x,u),
\qquad
\widetilde{u}=\psi(x,u)
\label{point_transform}
\end{equation}
has a long history. S. Lie showed, \cite{Lie1883}, that equations of the form
\begin{equation}
u_{xx}= A_3(x,u)\,u_x^3+A_2(x,u)\,u_x^2+A_1(x,u)\,u_x+A_0(x,u)
\label{LLT}
\end{equation}
generate an invariant subclass in the class (\ref{main}) with respect to the changes of variables
(\ref{point_transform}). He also showed that equation (\ref{LLT}) is linearizable by means of the
transformation (\ref{point_transform}) whenever the following system
\[
\fl
\quad
\begin{array}{lclclcl}
U_x &=& U\,V + A_0\,A_3 -\case{1}{3}\,A_{1,u} + \case{2}{3}\,A_{2,x},
&~~&
U_u &=& U^2 - A_2\,U + A_3\,V + A_{3,x} +A_1\,A_3,
\\
V_x &=& V^2 - A_0\,U + A_1\,V -A_{0,u} + A_0\,A_2,
&&
V_u &=& U\,V+ A_0\,A_3  +\case{1}{3}\,A_{2,u}-\case{2}{3}\,A_{1,u}
\end{array}
\]
is compatible. The compatibility condition of this system is equivalent to
the system $L_1\equiv 0$, $L_2 \equiv 0$ for the functions
\begin{eqnarray}
L_1 &=& 3\,A_{0,uu}-2\,A_{1,xu}+A_{2,xx}+3\,A_3 A_{0,x} - 3\,A_2 A_{0,u}
-3\,A_1 A_{1,u}
\nonumber
\\
&& - A_1 A_{2,x}
- 3\,A_0 A_{2,u}+6\,A_0 A_{3,x},
\label{L1}
\\
L_2 &=& A_{1,uu} -2 \,A_{2,xu} + 3\,A_{3,xx} -6\,A_3 A_{0,u} +3 \,A_3 A_{1,x}
+2\,A_2 A_{1,u}
\nonumber
\\
&&-2\, A_2 A_{2,x}
+ 3\,A_1 A_{3,x}-3\,A_0 A_{3,u}.
\label{L2}
\end{eqnarray}
These functions were found by R. Liouville in the first systematic study of the equivalence problem for
equations (\ref{LLT}) with respect to transformations (\ref{point_transform}), \cite{Liouville1889}.
Liouville found series of relative and absolute invariants and pointed out the procedure for generating
of invariants of higher orders.

A. Tresse used S. Lie's infinitesimal method to find differential invariants of equations (\ref{LLT}) and
(\ref{main}), \cite{Tresse1894,Tresse1896}. Papers \cite{Yumaguzhin2010} and \cite{Kruglikov2008} are devoted
to the modern exposition of Tresse's approach. Let us note one of the results of \cite{Tresse1894}: if either
$L_1$ or $L_2$ is not equal to zero, then there exists a change of variables (\ref{point_transform}) that maps
equation (\ref{LLT}) into equation with  $L_1 \not = 0$ and $L_2 \equiv 0$; moreover, two equations with
$L_2 \equiv 0$ are equivalent with respect to the pseudogroup (\ref{point_transform}) whenever they are
equivalent with respect to the pseudogroup
\begin{equation}
\widetilde{x} = \varphi(x),
\qquad
\widetilde{u} = \psi(x,u).
\label{fibre_preserving}
\end{equation}

{\'E}. Cartan developed the equivalence method \cite{Cartan_a}--\cite{Cartan_e} and apllied it to equation
(\ref{LLT}) in the paper \cite{Cartan_f}, however he studied there differnetial geometry of projective
connections, not equations  themselves, (see \cite[\S~8]{Cartan_f}).

The paper \cite{HsuKamran1989} applies Cartan's method to the equivalence problem for equations (\ref{main})
with respect to the pseudogroup (\ref{fibre_preserving}).

Subclass (\ref{LLT}) contains Painlev{\'e}'s equations \cite{Painleve1900,Painleve1902},  which are of interest
since they often appear in study of invariant solutions of untegrable nonlinear equations and in other
applications of physical importance \cite{Conte1999}. In \cite{KamranLambShadwick1985,HsuKamran1989,Kamran1989}
Cartan's method was used to find the necessary and sufficient conditions of equivalence of equation (\ref{LLT})
to the first or to the second Painlev{\'e} equation. In \cite{BabichBordag1999} it was
proven that all the Painlev{\'e} equations can be tranformed to the form
\begin{equation}
u_{xx} = A_0(x,u).
\label{generalizedEmdenFowler}
\end{equation}
In this paper it was also established that the point transformations that preserve
subclass (\ref{generalizedEmdenFowler}) have the form ~
$\widetilde{x}=\varphi(x)$, $\widetilde{u} = \sqrt{\left|\varphi^{\prime}(x)\right|}\,u+\chi(x)$.
In [2] R. Liouville's results were applied to study equations (\ref{LLT})
that are reductions of systems of differential equations describing chaotic dynamics in physics of the atmosphere and
and in chemical kynetics.

In the present paper we use Cartan's method \cite{Cartan_a}--\cite{Cartan_e},
\cite{Gardner1989,Kamran1989,Olver1995} to solve the equivalence problem for equations (\ref{main}) with respect
to transformations (\ref{point_transform}).

%%%%%%%%%%%%%%%%%%%%%%%%%%%%%%%%%%%%%%%%%%%%%%%%%%%%%%%%%%%%%%%%%%%%%%%%%%%%%%%%%%%%%%%%%%%%%%%%%%%%%%%%%%

\section{The solution of the equivalence problem in the case %\\
$F_{u_xu_xu_xu_x}\not=0$.}

All considerations in the paper are local, all maps are supposed to be real--analytic.

Equation  (\ref{main}) is a  submanifold in the bundle $J^2(\pi)$ of the second order jets of local sections of
the bundle $\pi \colon \mathbb{R}\times\mathbb{R} \rightarrow \mathbb{R}$,
$\pi \colon (x,u) \mapsto x$. Local coordinates on   $J^2(\pi)$ are $(x,u,u_x,u_{xx})$.
The pseudogroup of local diffeomorphisms
$\Phi \colon \mathbb{R}\times\mathbb{R} \rightarrow \mathbb{R}\times\mathbb{R}$,
$\Phi\colon (x,u) \mapsto (\widetilde{x},\widetilde{u})$,
acts on $\mathbb{R}\times\mathbb{R}$.
The second prolongation
$\Phi^{(2)} \colon J^2(\pi) \rightarrow J^2(\pi)$,
$\Phi^{(2)} \colon (x,u,u_x,u_{xx}) \mapsto
(\widetilde{x},\widetilde{u},\widetilde{u}_{\widetilde{x}},
\widetilde{u}_{\widetilde{x}\widetilde{x}})$
of a diffeomorphism $\Phi$ is defined as follows:
\[
\fl
(\Phi^{(2)})^{*}
\left(
\begin{array}{c}
d \widetilde{x}
\\
d \widetilde{u} - \widetilde{u}_{\widetilde{x}}\,d\widetilde{x}
\\
d \widetilde{u}_{\widetilde{x}}
 - \widetilde{u}_{\widetilde{x}{\widetilde{x}}}\,d\widetilde{x}
\end{array}
\right)
=
B
\left(
\begin{array}{c}
d x
\\
d u- u_x\,d x
\\
d u_x
 - u_{xx}\,d x
\end{array}
\right),
\quad
B\in \EuScript{H} = \left\{\left(
\begin{array}{ccc}
b_1 & b_2 & 0
\\
0 & b_3 & 0
\\
0 & b_5 & b_4
\end{array}
\right)
\in
\mathrm{GL}(3)
\right\}.
\]
Diffeomorphisms $\Phi^{(2)}$ constitute the pseudogroup
$\mathrm{Cont}_0(J^2(\pi))$ of point transformations of the bundle
$J^2(\pi)$.
When a superposition of two local diffeomorphisms from pseudogroup  $\mathrm{Cont}_0(J^2(\pi))$ is defined,
this superposition belongs to $\mathrm{Cont}_0(J^2(\pi))$ as well. Therefore the forms
\[
\left(
\begin{array}{c}
\Omega_1
\\
\Omega_2
\\
\Omega_3
\end{array}
\right)
=
B
\cdot
\left(
\begin{array}{c}
d x
\\
d u- u_x\,d x
\\
d u_x
 - u_{xx}\,d x
\end{array}
\right)
\]
are invariant with respect to the lifts
$\Psi \colon J^2(\pi) \times \EuScript{H} \rightarrow J^2(\pi) \times \EuScript{H}$
of diffeomorphisms from the pseudogroup  $\mathrm{Cont}_0(J^2(\pi))$.
Two equations (\ref{main}) are locally equivalent with respect to  $\mathrm{Cont}_0(J^2(\pi))$
whenever the restrictions
$\omega_i = \Omega_i \vert {}_{u_{xx}=F(x,u,u_x)}$ of forms  $\Omega_i$ onto these equations
are equivalent with respect to a diffeomorphism
$\Psi \colon J^2(\pi) \times \EuScript{H} \rightarrow J^2(\pi) \times \EuScript{H}$:   ~
$\Psi^{*}(\widetilde{\omega}_i) = \omega_i$,
$i \in \{1,2,3\}$. Hence we obtain $\EuScript{H}$-valued equivalence problem
for the collection of 1-forms
$\boldsymbol{\omega}$
$=\{\omega_1,\omega_2,\omega_3\}$
(see \cite[Def.~9.5]{Olver1995}).
In accordance with Cartan's method, to solve this problem we analyze the structure equations for forms
 $\omega_i$, that is, the expressions for the exterior differentials $d\omega_i$ via  $\boldsymbol{\omega}$.

The structute equation for form  $\omega_2$ is
\[
d \omega_2 = \eta \wedge \omega_2 + b_3b_1^{-1}b_4^{-1} \,\omega_1 \wedge \omega_3,
\]
where
$\eta = db_3\,b_3^{-1} - b_5 b_1^{-1}b_4^{-1}\,\omega_1 + r\,\omega_2
+b_2 b_1^{-1}b_4^{-1}\,\omega_3$, and  $r$ is an arbitrary constant.
Since forms  $\omega_i$ and their differentials are invsriant with respect to  $\Psi$, function
$b_3b_1^{-1}b_4^{-1}$ is invariant as well. We can normalize it, that is, to put it
equal to any non-zero constant, see \cite[Prop.~9.11]{Olver1995}. In the case $b_3b_1^{-1}b_4^{-1} = 1$ we get
$b_3=b_1 b_4$. Atfer this normalization the structure equations acquire the form
\begin{eqnarray*}
d \omega_1 &=& \eta_1 \wedge \omega_1 + \eta_2 \wedge \omega_2,
\\
d \omega_2 &=& \eta_3 \wedge \omega_2 + \omega_1 \wedge \omega_3,
\\
d \omega_3 &=& \eta_4 \wedge \omega_2 + (\eta_3-\eta_1) \wedge \omega_3,
\\
d \eta_1 &=& -2\,\eta_4 \wedge \omega_1 + \eta_5 \wedge \omega_2 - \eta_2 \wedge \omega_3,
\\
d \eta_2 &=& (\eta_1 -\eta_3) \wedge \eta_2 + \eta_5 \wedge \omega_1
+\case{1}{6}\,F_4\,b_1^{-1} b_4^{-3}\,\omega_2 \wedge \omega_3,
\end{eqnarray*}
where forms  $\eta_1$, ... , $\eta_4$
depend on differentials of the remaining non-normalized parameters of group $\EuScript{H}$,  form $\eta_5$
is obtained by means of the procedure of prolongation of the structure equations \cite[Ch.~12]{Olver1995},
and where we use the notation $F_k = \left(\frac{\partial}{\partial u_x}\right)^k F$, ~ $k \in \mathbb{N}$.

\vskip 5 pt
The further analysis divides on two cases:  case $\mathscr{A}$ corresponds to the condition $F_4\not\equiv 0$,
and case $\mathscr{B}$ corresponds to $F_4 \equiv0$.

\vskip 5 pt
In case  $\mathscr{A}$ we can shrink, if it is necessary, the domain of diffeomorphism $\Psi$, %and
therefore we can assume that $F_4\not=0$. Then normalization
$F_4b_1^{-1} b_4^{-3}=1$ yields $b_1 = F_4\,b_4^{-3}$.
After this we get
\begin{eqnarray*}
d \omega_1 &=& \eta_1 \wedge \omega_1 +\eta_2 \wedge \omega_2
+\case{1}{2}\,(5\,b_2b_4^3+F_5)b_4^{-1} F_4^{-1}\,\omega_1 \wedge \omega_3,
\\
d \omega_2 &=& \case{2}{3}\,\eta_1 \wedge \omega_2 + \omega_1 \wedge \omega_3
\end{eqnarray*}
with new forms $\eta_1$ and $\eta_2$. Then normalization $b_2 = - \frac{1}{5}\,F_5 \,b_4^{-3}$
gives the following structure equations
\begin{eqnarray*}
d \omega_1 &=& \eta_1 \wedge \omega_1
+\case{1}{25}\,(5\,F_4F_6-6\,F_5^2) F_4^{-2}b_4^{-2}\,\omega_2 \wedge \omega_3,
\\
d\omega_2 &=& \case{2}{3}\,\eta_1 \wedge \omega_2 +\omega_1\wedge \omega_3,
\\
d\omega_3 &=& \eta_2 \wedge \omega_2 - \case{1}{3}\,\eta_1 \wedge \omega_3
+\case{1}{2}\,b_4^2F_4^{-2}\,\left(F_4\,b_5+(D_x(F_4)+2\,F_1F_4)\,b_4 \right)\,
\omega_1 \wedge \omega_3,
\end{eqnarray*}
where we denote
$
D_x = \frac{\partial}{\partial x} + u_x\,\frac{\partial}{\partial u}
+ F\,\frac{\partial}{\partial u_x}
$.
We put the coefficient at $\omega_1 \wedge \omega_3$ in the third equation equal to zero
and obtain
$b_5 = - b_4 F_4^{-1}F\,(D_x(F_4)+2\,F_1F_4)$.

\vskip 5 pt

The further analysis divides on two subcases:
case $\mathscr{A}_1$ corresponds to the conditions
$5\,F_4F_6-6\,F_5^2 \not \equiv 0$,
while case $\mathscr{A}_2$ corresponds to the condition $5\,F_4F_6-6\,F_5^2\equiv 0$.

\vskip 5 pt

In the case $\mathscr{A}_1$ we can assume without loss of generality that
$5\,F_4F_6-6\,F_5^2 \not = 0$. Then we can put the coefficient at  $\omega_2 \wedge \omega_3$ equal to
$\case{1}{25}$. This yields $ b_4 = F_4^{-1}\,\sqrt{\vert 5\,F_4F_6 -6\,F_5^2 \vert}$.
After this normalization all the parameters of group $\EuScript{H}$ are defined as functions on $J^2(\pi)$,
and the structure equations acquire the form
\begin{eqnarray*}
d\omega_1
&=&
G_1\,\omega_{1}\wedge \omega_{2}
+ \left(G_2\,\omega_1 +\case{1}{25}\,\omega_2\right) \wedge \omega_3,
\\
d \omega_2
&=&
G_3\,\omega_1 \wedge \omega_2
+\left(\omega_1 + \case{2}{3} \,G_2\,\omega_2\right) \wedge \omega_3,
\\
d\omega_3
&=&
G_4\,\omega_1 \wedge \omega_2
+ \left(G_5\,\omega_2 - \case{1}{2}\,G_3\,\omega_1 \right)\wedge \omega_3,
\end{eqnarray*}
where  $G_1$, ... , $G_5$ are defined as
\begin{eqnarray}
\fl
G_1
&=&
\case{3}{10}\,
F_5 F_4^{-4} V_{1,x}
+
\case{3}{10}\,(u_x\,F_5 +5\,F_4)\,F_4^{-4}\,V_{1,u}
-\case{1}{5}\,V_1\,F_4^{-4}\,(F_{5,x}+u_x\,F_{5,u})
+
\case{6}{5}\,
G_2\,F\,F_5 F_4^{-5} V_1^{3/2}
\nonumber
\\
\fl &&
+\case{1}{5}\,V_1\,F_4^{-5}
\left(
\left(
5\,G_2\,V_1^{1/2} -3 \,F_5
\right)\,F_{4,x}
+
\left(
5\,(u_x\,G_2\,V_1^{1/2} - 4\,F_4) -3\,u_x\,F_5
\right)\,F_{4,u}
\right)
\nonumber
\\
\fl &&
+
2\,G_2\,F_1 V_1^{3/2}F_4^{-4}
-\case{2}{50}\,F\,V_1^2 F_4^{-5},
\nonumber
\\
\fl
G_2
&=&
\case{3}{10}\,V_1^{-3/2}\,
\left(
5\,F_4\,V_{1,u_x} - 14\,F_5\,V_1
\right),
\nonumber
\\
\fl
G_3
&=&
V_1^{1/2}F_4^{-5}\,
(
4\,V_1\,(F_{4,x}+u_x\,F_{4,u})
-F_4\,(V_{1,x} + u_x\,V_{1,u})
+\case{6}{5}\,F\,F_5\,V_1
+2\,F_1\,F_4\,V_1
-\case{2}{3}\,G_2\,F\,V_1^{3/2}
),
\nonumber
\\
\fl
G_4
&=&
-V_1^3\,F_4^{-9}\,
\left(
F_5\,(D_x(F)+3\,F\,F_1)
+2\,F_4\,(F_{5,x}+u_x\,F_{5,u})
+3\,F_1\,F_4\,F_{4,x}
\right.
\nonumber
\\
\fl
&&
\left.
+(3\,u_x\,F_1+F)\,F_{4,u}
+
2\,F_4\,D_x(F)
+u_x^2\,F_{4,uu}
-F_4\,F_u
+2\,F_1^2\,F_4
+F_{4,xx}
+2\,u_x\,F_{4,xu}
\right)
\nonumber
\\
\fl
&&-\case{1}{5}\,F^2\,V_1^3\,F_4^{-10}\,
(6\,F_5^2 +V_1)
\nonumber
\\
\fl
G_5
&=&
-5\,G_1
+\case{8}{5}\,
F_5\,F_4^{-4}\,(V_{1,x}+u_xV_{1,u})
+8\,F_4^{-3}\,V_{1,u}
+\case{4}{3}\,V_1\,F_4^{-5}\,F_{4,x}\,(4\,G_2\,V_1^{1/2}-3\,F_5)
\nonumber
\\
\fl
&&
-\case{4}{3}\,V_1\,F_4^{-5}\,F_{4,u}\,(u_x\,(3\,F_5-4\,G_2\,V_1^{1/2})+15\,F_4)
+\case{12}{25}\,F\,F_5^2\,V_1\,F_4^{-5}
\nonumber
\\
\fl
&&
+\case{8}{5}\,V_1\,F_5\,F_4^{-5}\,(F_1\,F_4+4\,G_2\,F\,V_1^{1/2})
+2\,F_2\,V_1\,F_4^{-3}
+\case{32}{3}\,G_2\,F_1\,V_1^{3/2}F_4^{-4}
\label{Invariants_AA}
\end{eqnarray}
with  $V_1=\left| 5\,F_4\,F_6-6\,F_5^2\right|$.
The functions $G_1$, ... , $G_5$ are invariants of equation (\ref{main}) with respect to
$\mathrm{Cont}_0(J^2(\pi))$. All the other invariants can be obtained from
$G_1$, ... , $G_5$ by applying the invariant derivatives
\begin{eqnarray}
\mathbb{D}_1
&=&
V_1^{3/2}\,F_4^{-4}\,D_x,
\nonumber
\\
\mathbb{D}_2
&=&
\case{1}{5}\,V_1\,\left(
F_5\,F_4^{-4}\,D_x + 5\,F_4^{-3}\,\case{\partial}{\partial u}
+5\,F_4^{-4}\,(D_x(F_4)+2\,F_1\,F_4)\case{\partial}{\partial u_x}
\right).
\label{InvDiff_AA}
\\
\mathbb{D}_3
&=&
F_4\,V_1^{-1/2}\,\case{\partial}{\partial u_x},
\nonumber
\end{eqnarray}
The operators $\mathbb{D}_1$ and $\mathbb{D}_2$ are defined by the requirement that
$dZ=\mathbb{D}_1(Z)\,\omega_1 + \mathbb{D}_2(Z)\,\omega_2 + \mathbb{D}_3(Z)\,\omega_3$
hold for an arbitrary function $Z(x,u,u_x)$.

The s${}^{th}$ order classifying manifold associated with forms $\boldsymbol{\omega}$,
in the case $\mathscr{A}_1$  has the form
\begin{equation}
\fl
\mathcal{C}^{(s)}_{\mathscr{A}_1}(\boldsymbol{\omega},\mathcal{U}) =
\left\{
\mathbb{D}_1^i\mathbb{D}_2^j\mathbb{D}_3^k(G_m(x,u,u_x))\,\,\vert\,\,
0\le i+j+k \le s, 1\le m\le 5, (x,u,u_x) \in \mathcal{U}
\right\},
\label{classifying_set_A_1}
\end{equation}
where  $\mathcal{U} \subset J^1(\pi)$ is an open subset such that
$5\,F_4F_6-6\,F_5^2 \not = 0$  everywhere in $\mathcal{U}$, and operators $\mathbb{D}_i$
are defined by equations (\ref{InvDiff_AA}).
Since all the invariants  (\ref{Invariants_AA}) depend on three variables $x$, $u$, $u_x$,
the number of functionally-independent invariants is no more than 3, and to formulate the solution
of the equivalence problem it is enough to consider 3${}^{rd}$ oredr classifying manifolds.

\vskip 5 pt

The analysis of the case $\mathscr{A}_2$ disparts on two subcases depending on whether one of the two
conditions $F_5 \equiv 0$ or $F_5 \not\equiv 0$ hold.
The first condition gives equations
\begin{equation}
\hspace{15pt}
u_{xx} = A_4(x,u)\,u_x^4+A_3(x,u)\,u_x^3+A_2(x,u)\,u_x^2+A_1(x,u)\,u_x+A_0(x,u),
\label{A4}
\end{equation}
with $A_4=F_4\not= 0$,
the second condition together with the defining identity $5\,F_4F_6-6\,F_5^2\equiv 0$ of subcase $\mathscr{A}_2$
gives  equations
\begin{equation}
\fl
u_{xx} = \left(B_1(x,u)\,u_x+B_0(x,u)\right)^{-1}+A_3(x,u)\,u_x^3+A_2(x,u)\,u_x^2+A_1(x,u)\,u_x+A_0(x,u)
\label{rational_1}
\end{equation}
with $B_1 \not = 0$.  However, the analysis of the second subcase can be reduced to the first one.

\vskip 7 pt

\noindent
{\sc Lemma   1}.
{\it
Each equation (\ref{rational_1}) can be mapped into equation (\ref{A4}) by means of a transformation from}
$\mathrm{Cont}_0(J^2(\pi))$.

\vskip 5 pt
\noindent
{\it Proof}:
consider 1-form $\mu = B_1\,du +B_0\,dx$. For its differential we have
\[
d\mu = (B_{1,x}-B_{0,u})\,dx \wedge du = (B_{1,x}-B_{0,u})\,B_1^{-1}\,dx\wedge \mu.
\]
This implies that form $\mu$ meets conditions of Frobenius' theorem
\cite[Th.~2.4.2]{Stormark2000}. Therefore there exists function $U=U(x,u)$ such that
$\mu \equiv 0 \,\,\mathrm{mod}\,\, d U$.
The direct check shows that the change of variables  $\widetilde{x} = U(x,u)$, $\widetilde{u} = x$
maps equation (\ref{rational_1}) into equation of the form (\ref{A4}).
\hfill {\sc qed}

\vskip 5 pt
\noindent
Moreover, without loss of generality we can put $A_3 \equiv 0$ in equation (\ref{A4}).

\vskip 5 pt
\noindent
{\sc Lemma 2}:
{\it
Each equation from the class  (\ref{A4}) can be mapped into equation from
the same class with} $A_3\equiv0$ {\it  via a transformation from} $\mathrm{Cont}_0(J^2(\pi))$.

\vskip 5 pt
\noindent
{\it Proof}:
consider  1-form $\nu = du -3\,A_3\,A_4^{-1}\,dx$.
If $(A_3A_4^{-1})_u = 0$, then $d\nu = 0$, otherwise
$d\nu = 3\,(A_3A_4^{-1})_u\,dx \wedge du = 3\,(A_3A_4^{-1})_u\,dx \wedge \nu$.
In both cases form $\nu$ meets the conditions of Frobenius' theorem,
hence there exists function $U=U(x,u)$ such that $\nu \equiv 0 \,\,\mathrm{mod}\,\, d U$.
The direct check shows that the change of variables  $\widetilde{x} = U(x,u)$, $\widetilde{u} = x$
maps equation from the class (\ref{A4}) into equation from the same class with  $A_3 \equiv 0$.
\hfill {\sc qed}

\vskip 7 pt

For equation
\begin{equation}
u_{xx} = A_4(x,u)\,u_x^4+A_2(x,u)\,u_x^2+A_1(x,u)\,u_x+A_0(x,u)
\label{A4_A3_0}
\end{equation}
after the above normalizations  we get the structure equations for the forms  $\omega_i$
\begin{eqnarray*}
d \omega_1 &=& \eta_1 \wedge \omega_1,
\\
d\omega_2 &=& \case{2}{3}\,\eta_1 \wedge \omega_2 +\omega_1\wedge \omega_3,
\\
d\omega_3 &=&  - \case{1}{3}\,\eta_1 \wedge \omega_3 + \dots \,\omega_1 \wedge \omega_2
+\case{1}{18}\,b_4^2\,u_x^2Z_0^{-2}\,
\omega_2 \wedge \omega_3,
\end{eqnarray*}
where $Z_0 = A_4u_x\,\left|18\,A_4^2u_x^2+3\,A_2A_4+A_{4,u} \right|^{-1/2}$.
We normalize the coefficient at  $\omega_2 \wedge \omega_3$ in the last equation by putting $b_4=Z_0\,u_x^{-1}$.
Then all the parameters of group $\EuScript{H}$ are defined as functions on  $J^2(\pi)$. We obtain
\begin{eqnarray}
\fl
\omega_1 &=& 24\,Z_0^{-3}\,A_4\,u_x^3\,dx,
\nonumber
\\
\fl
\omega_2 &=& 24\,Z_0^{-2}\,A_4\,u_x^2\,(du-u_x dx),
\label{old_omegas}
\\
\fl
\omega_3 &=& Z_0\,u_x^{-1}\,(du_x  - (A_4^2u_x^3+A_2A_4u_x+2\,A_1A_4+A_{4x})\,du
-(A_0A_4 - (A_1A_4+A_{4x})\,u_x)\,dx,
\nonumber
\end{eqnarray}
and
$d\omega_1 = 54\,Z_0\,\omega_1 \wedge \omega_2 + \dots \,\omega_1 \wedge \omega_3$,
that is, function $Z_0$ is an invariant of equations (\ref{A4_A3_0}) with respect to the
transformations from the pseudogroup $\mathrm{Cont}_0(J^2(\pi))$.
Since forms $\boldsymbol{\omega}$ and fucntion $Z_0$ are invariant, we can
take forms
$\hat{\omega}_1 = \case{1}{24}\,Z_0^3\,\omega_1$,
$\hat{\omega}_2 = \case{1}{24}\,Z_0^2\,\omega_2$ и
$\hat{\omega}_3 = Z_0^{-1}\,\omega_3$
instead of forms (\ref{old_omegas}).

In what follows we return to the previous notation, that is we will write  $\omega_i$ instead of $\hat{\omega}_i$.
The structure equations for the new forms are
\begin{eqnarray*}
%\fl
d\omega_1
&=&
-(3+Z_1\,u_x^{-2}+3\,Z_2\,u_x^{-3})\,\omega_1\wedge\omega_2 - 3\,\omega_1\wedge \omega_3,
\\
%\fl
d\omega_2 &=& (3+Z_1\,u_x^{-2}+2\,Z_2\,u_x^{-3}+2\,Z_3\,u_x^{-4})\,\omega_1\wedge\omega_2
+(\omega_1 -2\,\omega_2) \wedge\omega_3,
\\
%\fl
d\omega_3 &=& (Z_2 u_x^{-3}-3\,Z_3 u_x^{-4}+Z_5\,u_x^{-5}+Z_4u_x^{-6})\,\omega_1\wedge\omega_2
+(2-Z_2u_x^{-3}-Z_3u_x^{-4})\,\omega_1\wedge\omega_3
\\
%\fl
&&
+(2-Z_2u_x^{-3})\,\omega_2\wedge\omega_3,
\end{eqnarray*}
where
\begin{eqnarray*}
Z_1 &=& A_{4,u}A_4^{-1}+3\,A_2,
\\
Z_2 &=& A_{4,x}A_4^{-2}+2\,A_1A_4^{-1},
\\
Z_3 &=& A_0A_4^{-1},
\\
Z_4 &=& (Z_{3,u}-Z_{2,x}+3\,A_1Z_2-A_2Z_3)\,A_4^{-1}+Z_3\,A_{4,u}A_4^{-2}-2\,Z_2^2,
\\
Z_5 &=& (A_{1,u}-A_{2,x})\,A_4^{-2}-(Z_{2,u}-3\,A_2Z_2)\,A_4-Z_1Z_2.
\end{eqnarray*}

\vskip 5 pt

The further analysis of the class of equations (\ref{A4_A3_0})
divides onto five cases in dependence on whether the functions  $Z_1$, ..., $Z_5$ are equal to zero.

\vskip 3 pt
%\noindent
The case $\mathscr{A}_{21}$ corresponds to condition  $Z_1 \not =0$.
In this case the structure equations acquire the form
\begin{eqnarray*}
\fl
d\omega_1 &=& -(3\,H_1P_1^3+P_1^2+3)\,\omega_1\wedge\omega_2-3\,\omega_1\wedge\omega_3,
\\
\fl
d\omega_2 &=& (2\,H_2P_1^4+2\,H_1P_1^3+P_1^2+3)\,\omega_1\wedge\omega_2
+(\omega_1-2\,\omega_2) \wedge\omega_3,
\\
\fl
d\omega_3 &=& P_1^3\,(H_3P_1^3+H_4P_1^2-3\,H_2P_1+H_1)\,\omega_1\wedge\omega_2
-(H_2P_1^4+H_1P_1^3-2)\,\omega_1\wedge\omega_3
-(H_1P_1^3-2)\,\omega_2\wedge\omega_3,
\end{eqnarray*}
where the unique invariant that depends on $x$, $u$, and $u_x$ is defined as
$P_1 = Z_1^{1/2}\,u_x^{-1}$, and its differential is
\[
\fl
dP_1 = -\case{1}{2}\,P_1\,(2\,H_2P_1^4-H_6P_1^3-H_5P_1^2+2)\,\omega_1
-\case{1}{2}\,P_1\,(2\,H_1P_1^3-H_5P_1^2+2)\,\omega_2
-P_1\,\omega_3.
\]
All the invariants that depend on $x$ and  $u$ only are
\begin{eqnarray*}
H_1&=&Z_2 Z_1^{-3/2},
\\
H_2 &=& Z_3 Z_1^{-2},
\\
H_3 &=& Z_4 Z_1^{-3},
\\
H_4 &=& Z_5 Z_1^{-5/2},
\\
H_5 &=& (Z_{1,u}-2\,A_2Z_1)\,Z_1^{-2}A_4^{-1},
\\
H_6 &=& (Z_{1,x}-2\,A_1Z_1)\,Z_1^{-5/2}A_4^{-1}.
\end{eqnarray*}
The invariant derivatives
$
\mathbb{D}_1 = A_4^{-1}Z_1^{-3/2}\,\case{\partial}{\partial x}
$
and
$
\mathbb{D}_2 = A_4^{-1}Z_1^{-1}\,\case{\partial}{\partial u}
$
are defined by the requirement that equation
$
dY=P_1^2\,\left((P_1\mathbb{D}_1(Y)+\mathbb{D}_2(Y))\,\omega_1
+\mathbb{D}_1(Y)\,\omega_2\right)
$
holds for every function  $Y(x,u)$.
In the case $\mathscr{A}_{21}$  the second order classifying manifold associated with forms
$\boldsymbol{\omega}$ has the form
\begin{equation}
\mathcal{C}^{(2)}_{\mathscr{A}_{21}}(\boldsymbol{\omega},\mathcal{V}) =
\left\{
\mathbb{D}_1^i\mathbb{D}_2^j(H_m(x,u))\,\,\vert\,\,
0\le i+j \le 2, \,\, 1\le m\le 6, \,\, (x,u) \in \mathcal{V}
\right\},
\label{classifying_set_A_21}
\end{equation}
where $\mathcal{V} \subset J^0(\pi)$ is an open subset such that
$F_4\not = 0$, $F_5 \equiv 0$  and $Z_1\not = 0$   everywhere on it.

\vskip 5 pt
%\noindent
The case  $\mathscr{A}_{22}$ is defined by   $Z_1\equiv 0$ and $Z_2\not=0$. In this case the structure equations
\begin{eqnarray*}
\fl
d\omega_1 &=& 3\,\left((P_2^2+1)\,\omega_2+\omega_3\right)\wedge \omega_1,
\\
\fl
d\omega_2 &=& (2\,I_1P_2^4+2\,P_2^3+3)\,\omega_1\wedge\omega_2
+(\omega_1-2\,\omega_2)\wedge\omega_3,
\\
\fl
d\omega_3 &=& P_2^3\,(I_2P_2^3+I_3P_2^2-3\,I_1P_2+1)\,\omega_1\wedge\omega_2
-\left((I_1P_2^4+P_2^3-2)\,\omega_1+(P_2^3-2)\,\omega_2\right) \wedge\omega_3
\end{eqnarray*}
contain invariant $P_2 = Z_2^{1/3}\,u_x^{-1}$ that depends on $x$, $u$, and $u_x$,
with the differential
\[
\fl
dP_2 = -\case{1}{3}\,P_2\,(3\,I_1P_2^4-I_5P_2^3-I_4P_2^2+3)\,\omega_1
-\case{1}{3}\,P_3\,(3\,P_2^3+(I_3-I_4)P_2^2+3)\,\omega_2
-P_2\,\omega_3.
\]
Invariants that depend on  $x$ and $u$ only have the form
\begin{eqnarray*}
I_1 &=& Z_3 Z_2^{-4/3},
\\
I_2 &=& Z_4 Z_2^{-2},
\\
I_3 &=& Z_5 Z_2^{-5/3},
\\
I_4 &=& (A_{1,u}-A_{2,x})\,Z_2^{-5/3}A_4^{-2},
\\
I_5 &=& (Z_{2,x}-3\,A_1Z_2)\,Z_2^{-2}A_4^{-1}.
\end{eqnarray*}
The invariant derivatives \,
$\mathbb{D}_1 = A_4^{-1}Z_2^{-1}\,\case{\partial}{\partial x}$
and
$\mathbb{D}_2 = A_4^{-1}Z_2^{-2/3}\,\case{\partial}{\partial u}$
are defined by the requirement that equation
$
dY = P_2^2\, \left((P_2\mathbb{D}_1(Y) + \mathbb{D}_2(Y))\,\omega_1
  +  \mathbb{D}_1(Y)\,\omega_2\right)
$
holds for every function $Y(x,u)$.
In the case $\mathscr{A}_{22}$  the second order classifying manifold associated with forms  $\boldsymbol{\omega}$
has the form
\begin{equation}
\mathcal{C}^{(2)}_{\mathscr{A}_{22}}(\boldsymbol{\omega},\mathcal{V}) =
\left\{
\mathbb{D}_1^i\mathbb{D}_2^j(I_m(x,u))\,\,\vert\,\,
0\le i+j \le 2, \,\, 1\le m\le 5, \,\, (x,u) \in \mathcal{V}
\right\},
\label{classifying_set_A_22}
\end{equation}
where  $\mathcal{V} \subset J^0(\pi)$ is an open subset such that
$F_4\not = 0$, $F_5 \equiv 0$,  $Z_1\equiv 0$, and $Z_2\not = 0$ hold everywhere in it.

\vskip 5 pt
%\noindent
The subcase $\mathscr{A}_{23}$ is defined by conditions  $Z_1\equiv 0$, $Z_2\equiv 0$, $Z_3\not=0$.
In this case the structure equations have the form
\begin{eqnarray*}
d\omega_1 &=& 3\,(\omega_2+\omega_3)\wedge\omega_1,
\\
d\omega_2 &=& (2\,P_3^4+3)\,\omega_1\wedge\omega_2
+(\omega_1-2\,\omega_2)\wedge\omega_3,
\\
d\omega_3 &=& P_3^4\,(J_1P_3^2+J_2P_3-3)\,\omega_1\wedge\omega_2
-\left((P_3^4-2)\,\omega_1-2\,\omega_2\right)\wedge\omega_3,
\end{eqnarray*}
where $P_3 = Z_3^{1/4}\,u_x^{-1}$ and
\[
dP_3 = -\case{1}{4}\,P_3\,(4\,P_3^4-J_3P_3^3-J_1P_3^2+4)\,\omega_1
+\case{1}{4}\,P_3\,(J_1\,P_3^3-4)\,\omega_2
-P_3\,\omega_3.
\]
In this case the invariants that depend on $x$ and  $u$ are
\begin{eqnarray*}
J_1 &=& Z_4 Z_3^{-3/2},
\\
J_2 &=& Z_5 Z_3^{-5/4},
\\
J_3 &=& (A_4 Z_{3,x}+2\,Z_3A_{4,x})\,Z_3^{-7/4}A_4^{-2}.
\end{eqnarray*}
The invariant derivatives
$\mathbb{D}_1 = A_4^{-1}Z_3^{-3/4}\,\case{\partial}{\partial x}$
and
$\mathbb{D}_2 = A_4^{-1}Z_3^{-1/2}\,\case{\partial}{\partial u}$
are defined by equation
$
dY=P_3^2\,\left((P_3\mathbb{D}_1(Y)+\mathbb{D}_2(Y))\,\omega_1+\mathbb{D}_1(Y)\,\omega_2\right)
$.
In the case $\mathscr{A}_{23}$  the second order classifying manifold associated with forms  $\boldsymbol{\omega}$
has the form
\begin{equation}
\mathcal{C}^{(2)}_{\mathscr{A}_{23}}(\boldsymbol{\omega},\mathcal{V}) =
\left\{
\mathbb{D}_1^i\mathbb{D}_2^j(J_m(x,u))\,\,\vert\,\,
0\le i+j \le 2, \,\, 1\le m\le 3, \,\, (x,u) \in \mathcal{V}
\right\},
\label{classifying_set_A_23}
\end{equation}
where $\mathcal{V} \subset J^0(\pi)$ is an open subset such that
$F_4\not = 0$, $F_5 \equiv 0$,  $Z_1\equiv 0$, $Z_2\equiv 0$, and $Z_3\not = 0$  in all its points.

\vskip 5 pt
%\noindent
The subcase $\mathscr{A}_{24}$ is defined by conditions
$Z_1\equiv 0$, $Z_2\equiv 0$, $Z_3\equiv 0$, which imply $Z_4\equiv 0$, but $Z_5\neq 0$.
In this case the structure equations have the form
\begin{eqnarray*}
d\omega_1 &=& 3\,(\omega_2+\omega_3)\wedge\omega_1,
\\
d\omega_2 &=& 3\,\omega_1\wedge\omega_2
+(\omega_1-2\,\omega_2)\wedge\omega_3,
\\
d\omega_3 &=& -P_4^5\,\omega_1\wedge\omega_2
+2\,(\omega_1+\omega_2)\wedge\omega_3
\end{eqnarray*}
They contain only one invariant $P_4 = Z_5^{1/5}\,u_x^{-1}$ that depends on  $x$, $u$, and $u_x$,
and has differential
\[
dP_4 = \case{1}{5}\,P_4\,(K_2\,P_4^4+K_1P_4^2-5)\,\omega_1
+\case{1}{5}\,P_4\,(K_1\,P_4^2-5)\,\omega_2
-P_4\,\omega_3,
\]
which contains invariants
\begin{eqnarray*}
K_1 &=& (3\,A_4 Z_{5,u}+5\,Z_5A_{4,u})\,Z_5^{-7/5}A_4^{-2},
\\
K_2 &=& (3\,A_4 Z_{5,x}+5\,Z_5A_{4,x})\,Z_5^{-8/5}A_4^{-2}
\end{eqnarray*}
that depend on $x$ and $u$. The invariant derivatives
$\mathbb{D}_1 = A_4^{-1}Z_5^{-3/5}\,\case{\partial}{\partial x}$
and
$\mathbb{D}_2 = A_4^{-1}Z_5^{-2/5}\,\case{\partial}{\partial u}$
are defined by
\[
dY=P_4^2\,\left((P_4\mathbb{D}_1(Y)+\mathbb{D}_2(Y))\,\omega_1+\mathbb{D}_1(Y)\,\omega_2\right).
\]
In the case $\mathscr{A}_{24}$  the second order classifying manifold associated with forms  $\boldsymbol{\omega}$
has the form
\begin{equation}
\mathcal{C}^{(2)}_{\mathscr{A}_{24}}(\boldsymbol{\omega},\mathcal{V}) =
\left\{
\mathbb{D}_1^i\mathbb{D}_2^j(K_m(x,u))\,\,\vert\,\,
0\le i+j \le 2, \,\, 1\le m\le 2, \,\, (x,u) \in \mathcal{V}
\right\},
\label{classifying_set_A_24}
\end{equation}
where $\mathcal{V} \subset J^0(\pi)$ is an open subset such that
$F_4\not = 0$, $F_5 \equiv 0$,  $Z_1\equiv 0$, $Z_2\equiv 0$, $Z_3\equiv 0$, $Z_4 \equiv 0$, and $Z_5\not = 0$
everywhere in it.

\vskip 5 pt
%\noindent
Finally, the case  $\mathscr{A}_{25}$ is defined by the requirements
$Z_1 \equiv Z_2 \equiv Z_3 \equiv Z_4 \equiv Z_5 \equiv 0$ . In this case the structure equations
have the form
\begin{eqnarray*}
d\omega_1 &=& 3\,(\omega_2+\omega_3)\wedge\omega_1,
\\
d\omega_2 &=& 3\,\omega_1\wedge\omega_2
+(\omega_1-2\,\omega_2)\wedge\omega_3,
\\
d\omega_3 &=& 2\,(\omega_1 + \omega_2)\wedge\omega_3.
\end{eqnarray*}
Their coefficients are constant.  The same structure equations for forms $\boldsymbol{\omega}$
has the symmetry pseudogroup of equation
\begin{equation}
u_{xx}= u_x^4.
\label{A_25_normal_form}
\end{equation}

Therefore the class of equations (\ref{main}) such that  $F_{u_xu_xu_xu_x}\not=0$ is divided on invariant
subclasses $\mathscr{A}_{1}$, $\mathscr{A}_{21}$, ... , $\mathscr{A}_{25}$, and we construct invariants
coframes $\boldsymbol{\omega}$ for each subclass. Hence the solution of the equivalence problem for equations
(\ref{main}) is reduced to the restricted equivalence problem for invarinat coframes, see
\cite[Ch.~8, Ch.~14]{Olver1995}. The results of the above computations together with Theorem 14.24 from
\cite{Olver1995} give the following theorem:

\vskip 7 pt
\noindent
{\sc Theorem 1}:
{\it
Each equation (\ref{main}) such taht $F_{u_xu_xu_xu_x}\not=0$ can be mapped by means of a diffeomorphism from
the pseudogroup of point transformations (\ref{point_transform}) into an equation from one of the
invariant subclasses $\mathscr{A}_{1}$, $\mathscr{A}_{21}$, ... , $\mathscr{A}_{25}$.

The invariant subclass  $\mathscr{A}_{1}$ contains equations (\ref{main}) such that
\[
5\,F_{u_xu_xu_xu_x}F_{u_xu_xu_xu_xu_xu_x}-6\,F_{u_xu_xu_xu_xu_x}^2 \not = 0.
\]
Two equations from the invariant subclass $\mathscr{A}_{1}$
are locally equivalent whenever their classifying
manifolds (\ref{classifying_set_A_1}) are locally congruent.

Equations (\ref{main}) with $F_{u_xu_xu_xu_x}\not=0$, $5\,F_{u_xu_xu_xu_x}F_{u_xu_xu_xu_xu_xu_x}-6\,F_{u_xu_xu_xu_xu_x}^2 \equiv 0$
can be mapped into an equation of the form (\ref{A4_A3_0}).

The invariant subclass  $\mathscr{A}_{21}$ contains equations (\ref{A4_A3_0}) such that $Z_1\not=0$,
the invariant subclass $\mathscr{A}_{22}$ contains equations (\ref{A4_A3_0}) such that
    $Z_1\equiv 0$, $Z_2\not=0$,
the invariant subclass $\mathscr{A}_{23}$ contains equations (\ref{A4_A3_0}) such that
    $Z_1\equiv Z_2 \equiv 0$, $Z_3\not=0$,
the invariant subclass $\mathscr{A}_{24}$ contains equations (\ref{A4_A3_0}) such that
     $Z_1\equiv Z_2 \equiv Z_3 \equiv Z_4 \equiv 0$, $Z_5\not=0$,
the invariant subclass $\mathscr{A}_{25}$ contains equations (\ref{A4_A3_0}) such that
     $Z_1\equiv Z_2 \equiv Z_3 \equiv Z_4 \equiv Z_5 \equiv 0$.

Equations from the invariant subclasses  $\mathscr{A}_{21}$, ... , $\mathscr{A}_{24}$
are locally equivalent with respect to the pseudogroup of point transformations (\ref{point_transform})
whenever their classifying manifolds
(\ref{classifying_set_A_21}), (\ref{classifying_set_A_22}), (\ref{classifying_set_A_23}), and
(\ref{classifying_set_A_24}) are locally congruent.
Equations from the invariant subclass $\mathscr{A}_{25}$ are locally equivalent to equation (\ref{A_25_normal_form}).
}

%%%%%%%%%%%%%%%%%%%%%%%%%%%%%%%%%%%%%%%%%%%%%%%%%%%%%%%%%%%%%%%%%%%%%%%%%%%%%%%%%%%%%%%%%%%%%%%%%%%%%%%%%%%%%

\section{The solution of the equivalence problem in the case  $F_{u_xu_xu_xu_x}=0$.}

After the normalization $b_3=b_1 b_4$
and applying the procedure of prolongation
we get the structure equations of forms $\omega_i$
for equation (\ref{LLT}):
\begin{eqnarray}
d\omega_1 &=& \eta_1 \wedge \omega_1 + \eta_2 \wedge \omega_2,
\nonumber
\\
d\omega_2 &=& \eta_3 \wedge \omega_2 + \omega_1 \wedge \omega_3,
\nonumber
\\
d\omega_3 &=& (\eta_3-\eta_1) \wedge \omega_3 + \eta_4 \wedge \omega_2,
\nonumber
\\
d\eta_1 &=& -2\,\eta_4 \wedge \omega_1 + \eta_5 \wedge \omega_2- \eta_2 \wedge \omega_3,
\nonumber
\\
d\eta_2 &=& \eta_5 \wedge \omega_1 + (\eta_1-\eta_3) \wedge \eta_2,
\label{se_B}
\\
d\eta_3 &=& -\eta_4 \wedge \omega_1 + 2\,\eta_5 \wedge \omega_2+ \eta_2 \wedge \omega_3,
\nonumber
\\
d\eta_4 &=& \eta_4 \wedge \eta_1 + \eta_5 \wedge \omega_3
- \case{1}{3}\,b_1^{-3}b_4^{-1}\,(L_1+L_2\,u_x)\, \omega_1 \wedge \omega_2,
\nonumber
\end{eqnarray}
where
\[
\fl
L_1 =
3\,(A_{0,uu}+A_3 A_{0,x}-A_2 A_{0,u})
-2\,A_{1,xu}+A_{2,xx}
-A_1\,(3\,A_{1,u} + A_{2,x})
- 3\,A_0 (A_{2,u}-2\,A_{3,x}),
\label{L1}
\]
\[
\fl
L_2 = A_{1,uu}
-2 \,A_{2,xu}
+ 3\,(A_{3,xx} -A_3\,( A_{1,x}-2\,A_{0,u})+ A_1 A_{3,x}-A_0 A_{3,u})
+2\,A_2 (A_{1,u} -A_{2,x}).
\]
The further analysis depends on whether condition $L_1 \equiv L_2 \equiv 0$ holds.

\vskip 5 pt
In the case $\mathscr{B}_1$, when $L_1 \equiv L_2 \equiv 0 $,
the structure equations after the prolongation acquire the form
\begin{eqnarray*}
d\omega_1 &=& \eta_1 \wedge \omega_1 + \eta_2 \wedge \omega_2,
\\
d\omega_2 &=& \eta_3 \wedge \omega_2 + \omega_1 \wedge \omega_3,
\\
d\omega_3 &=& (\eta_3-\eta_1) \wedge \omega_3 + \eta_4 \wedge \omega_2,
\\
d\eta_1 &=& -2\,\eta_4 \wedge \omega_1 + \eta_5 \wedge \omega_2- \eta_2 \wedge \omega_3,
\\
d\eta_2 &=& \eta_5 \wedge \omega_1 + (\eta_1-\eta_3) \wedge \eta_2,
\\
d\eta_3 &=& -\eta_4 \wedge \omega_1 + 2\,\eta_5 \wedge \omega_2+ \eta_2 \wedge \omega_3,
\\
d\eta_4 &=& \eta_4 \wedge \eta_1 + \eta_5 \wedge \omega_3,
\\
d\eta_5 &=& \eta_2 \wedge \eta_4 + \eta_5 \wedge \eta_3.
\end{eqnarray*}
The direct check shows that any second-order linear ODE $u_{xx} = a_2(x)\,u_x+a_1(x)\,u+a_0(x)$, in particular,
equation $u_{xx}= 0$, has the same structure equations. Hence equations (\ref{LLT}) with
$L_1 \equiv L_2 \equiv 0$ are equivalent to equation $u_{xx}=0$, in agreement with S. Lie's and R. Liouville's
results.

\vskip 5 pt

We now turn to the analysis of the case $\mathscr{B}_2$, in which one of the functions
$L_1$  or $L_2$ is not equal to zero.
We put the coefficient at $\omega_1\wedge \omega_2$ in the last equation of (\ref{se_B})  equal to $\case{1}{3}$.
This yields $b_4 = b_1^{-3}\,(L_1+L_2\,u_x)$. Then we have
\begin{eqnarray*}
d\omega_1 &=& \eta_1 \wedge \omega_1 + \eta_2 \wedge \omega_2,
\\
d\omega_2 &=& -2\,\eta_2 \wedge \omega_2 + \omega_1 \wedge \omega_2
+ b_1^2\,(L_1+L_2\,u_x)^{-2}\,\left(b_2\,(L_1+L_2\,u_x)-b_1\,L_2\right)\,\omega_2\wedge\omega_3.
\end{eqnarray*}
Further, we put the coefficient at $\omega_2\wedge\omega_3$ in the last equation equal to zero and obtain
$b_2 = b_1\,L_2\,(L_1+L_2\,u_x)^{-1}$. Now we have $\omega_1 = b_1\,(L_1+L_2\,u_x)^{-1}\,(L_1\,dx+L_2\,du)$ and
$d\omega_1 = \eta_1 \wedge \omega_1$. Therefore $\omega_1$ meets the conditions of Frobenius' theorem, so there
are functions $X(x,u)$ and $\widetilde{b}_1(x,u,u_x,b_1,b_5)$ such that $\omega_1 = \widetilde{b}_1\, d X$.
Let   $U(x,u)$ be any function such that $dX \wedge dU \not =0$. Then we have
$\omega_1 = \widetilde{b}_1 d\widetilde{x}$ after the change of variables $\widetilde{x}= X(x,u)$,
$\widetilde{u}=U(x,u)$.
In the new variables, equality
$\omega_1 =
\widetilde{b}_1\,(\widetilde{L}_1+\widetilde{L}_2\,\widetilde{u}_{\widetilde{x}})^{-1}
\,(\widetilde{L}_1\,d\widetilde{x}+\widetilde{L}_2\,d\widetilde{u})$
should hold. This implies $\widetilde{L}_2=0$,
that is, we get Tresse's result, \cite{Tresse1896}:
for a nonlinearizable equation  (\ref{LLT})
there is a change of variables  (\ref{point_transform})
such that $\widetilde{L}_1 \not = 0$ and $\widetilde{L}_2\equiv 0$ in the new coordinates  $\widetilde{x}$,
$\widetilde{u}$. If it is necessary, we make this change of variables, therefore without loss of generality we
assume that conditions $L_1 \not = 0$, $L_2\equiv 0$ hold for equation (\ref{LLT}).

Now the structure equations have the form
\begin{eqnarray*}
\fl
d\omega_1 &=& \eta_1 \wedge \omega_1,
\\
\fl
d\omega_2 &=& - 2\,\eta_1 \wedge \omega_2 +\omega_1\wedge\omega_3,
\\
\fl
d\omega_3 &=& \eta_2 \wedge \omega_2
+\case{1}{3}\,b_1^{-1}L_1^{-1}\,
\left(
5\,b_1^3b_5 +6\,A_3L_1u_x^2+(4\,A_2L_1-L_{1,u})\,u_x
+2\,A_1L_1-L_{1,x}
\right)\,\omega_1 \wedge \omega_3
\\
\fl
&&
- 3\,\eta_1 \wedge \omega_3.
\end{eqnarray*}
We can assume
$b_5 =  - \case{1}{5}\,b_1^{-3}(6\,A_3L_1u_x^2+(4\,A_2L_1-L_{1,u})\,u_x +2\,A_1L_1-L_{1,x})$.
Then we obtain
\begin{eqnarray*}
d\omega_1 &=& \eta_1 \wedge \omega_1,
\\
d\omega_2 &=& - 2\,\eta_1 \wedge \omega_2 +\omega_1\wedge\omega_3,
\\
d\omega_3 &=& - 3\,\eta_1 \wedge \omega_3
+\case{4}{5}\,b_1^2L_1^{-2}\,
\left(
3\,A_3L_1u_x+A_2L_1 + L_{1,u}
\right)\,\omega_2 \wedge \omega_3
+\dots\,\omega_1 \wedge \omega_2.
\end{eqnarray*}

\vskip 5 pt

The further analysis seperates on two cases: the case $\mathscr{B}_{21}$ such that
$A_3 \not =0$, and the case $\mathscr{B}_{22}$ such that  $A_3 \equiv 0$.

\vskip 5 pt

In the case $\mathscr{B}_{21}$ we put
\begin{equation}
b_1 = L_1\,\left|3\,A_3L_1u_x+A_2L_1 + L_{1,u}\right|^{-1/2}.
\label{second_normalization}
\end{equation}
This yields
$
d\omega_1 = \case{3}{2}\,M_0^5\,\omega_1\wedge\omega_3 + \dots \, \omega_1\wedge\omega_2,
$
where $M_0 = A_3^{1/5}L_1^{3/5}\left|3\,A_3L_1u_x+A_2L_1 + L_{1,u}\right|^{-1/2}$,
and
\begin{eqnarray}
\fl
\omega_1 &=& M_0\,A_3^{-1/3}L_1^{2/5}\,dx,
\nonumber
\\
\fl
\omega_2 &=& M_0^{-2}A_3^{2/5}L_1^{1/5}\,(du-u_x dx),
\nonumber
\\
\fl
\omega_3 &=& M_0^{-3}A_3^{3/5}L_1^{-1/5}\,
\left(
d u_x
-\case{1}{5}\,\left(6\,A_3u_x^2+(4\,A_2-L_{1,u}L_1^{-1})\,u_x+2A_1-L_{1,x}L_1^{-1}\right) \, du
\right.
\nonumber
\\
\fl
&&
\left.
+ \case{1}{5}\,(A_3u_x^3-(A_2+L_{1,u}L_1^{-1}u_x^2-(3\,A_1+L_{1,x}L_1^{-1})\,u_x-5\,A_0)\,dx
\right).
\label{omegas_BBA}
\end{eqnarray}
Since forms $\boldsymbol{\omega}$ and function $M_0$ are invariant with respect to a diffeomorphism
$\Psi$, we can multiply, without loss of generality, the right hand sides of forms $\omega_1$, $\omega_2$, and $\omega_3$ by $M_0^{-1}$,
$M_0^2$, and  $M_0^3$, respectively. We denote the obtained forms as $\omega_1$, $\omega_2$, and $\omega_3$ again.
These forms have the following structure equations:
\begin{eqnarray*}
d\omega_1 &=& M_1 \,\omega_1 \wedge \omega_2,
\\
d\omega_2 &=& \case{2}{15}\,(P_5^2 +5\,M_1\,P_5+M_2)\,\omega_1 \wedge \omega_2
+\omega_1\wedge \omega_3,
\\
d\omega_3 &=& -\case{1}{675}\,\left(
2\,P_5^4+25\,M_1P_5^3+3\,(10\,M_2-M_3)\,P_5^2-M_4P_5-M_5
\right)\,\omega_1\wedge\omega_2
\\
&&
+\case{1}{5}\,(P_5^2+5\,M_1P_5+M_2)\,\omega_1\wedge\omega_3
+\case{1}{5}\,(4\,P_5+15\,M_1)\,\omega_2\wedge\omega_3,
\end{eqnarray*}
where
$P_5 = A_3^{-2/5}L_1^{-6/5}\,(3\,A_3L_1u_x+A_2L_1+L_{1,u})$ is the only differential invariant that depends
on $x$, $u$, and $u_x$, while
\[
\fl
dP_5 =\case{1}{90}\,(
10\,P_5^3+90\,M_1P_5^2+6\,(3\,M_2+M_3)\,P_5+15\,M_6
)\,\omega_1
+3\,\omega_3
+\case{1}{5}\,(2\,P_5^2+15\,M_1P_5+M_3)\,\omega_2,
\]
and invariants $M_1$, ... , $M_6$ depend on $x$ and $u$:
\begin{eqnarray*}
\fl
M_1
&=&
\case{1}{5}\,A_3^{-7/5}L_1^{-6/5}\,(L_1A_{3,u}-2\,A_3L_{1,u}),
\\
\fl
M_2
&=&
- A_3^{-4/5}L_1^{12/5}\,\left(
L_{1,u}^2+(2\,A_2+5\,M_1A_3^{2/5}L_1^{1/5})\,L_1L_{1,u}
+L_1^2(A_2^2-3\,A_1A_3-3\,A_{3,x}
\right.
\\
\fl &&\left.
+5\,M_1A_2A_3^{2/5}L_1^{1/5})
\right),
\\
\fl
M_3 &=&  A_3^{-4/5}L_1^{-7/5}\,\left(
5\,L_{1,uu} -12\,L_1^{-1}L_{1,u}^2
-(9\,A_2+25\,M_1A_3^{2/5}L_1^{1/5})\,L_{1,u}
-3\,A_3L_{1,x}
\right.
\\
\fl
&&
\left.
-L_1\,(
2\,A_2^2-6\,A_1A_3-5\,A_{2,u}+25\,M_1A_2A_3^{2/5}L_1^{1/5}
)
\right),
\\
\fl
M_4
&=&
A_3^{-1/5}L_1^{-8/5}\,\left(
90\,L_{1,xu}
+8\,A_3^{-1}L_1^{-2}L_{1,u}^3
+3\,A_3^{-1}L_1^{-1}(8\,A_2+25\,M_1A_3^{2/5}L_1^{1/5})\,L_{1,u}^2
\right.
\\
\fl
&&
\left.
-18\,L_1^{-1}\,(6\,L_{1,u}+A_2L_1)\,L_{1,x}
+135\,L_1A_{1,u}-180\,L_1A_{2,x}
+6\,A_3^{-1}\,(
4\,A_2^2+6\,A_1A_3
\right.
\\
\fl
&&
\left.
+25\,M_1A_2A_3^{2/5}L_1^{1/5}
+(10\,M_2-M_3)\,A_3^{4/5}L_1^{2/5}
)\,L_{1,u}
+A_3^{-1}L_1\,(8\,A_2^3+36\,A_1A_2A_3
\right.
\\
\fl
&&
\left.
-540\,A_0A_3^2
+6\,(10\,M_2-M_3)\,A_2A_3^{4/5}L_1^{4/5}
+75\,M_1A_2^2A_3^{2/5}L_1^{1/5}
)
\right),
\\
\fl
M_5
&=&
A_3^{2/5}L_1^{-9/5}\,\left(
135\,L_{1,xx}+2\,A_3^{-2}L_1^{-3}L_{1,u}^4
+A_3^{-2}L_1^{-3}\,(8\,A_2+25\,M_1A_3^{2/5}L_1^{1/5})\,L_{1,u}^3
-27\,A_1L_{1,x}
\right.
\\
\fl
&&
\left.
+3\,A_3^{-2}L_1^{-1}\,(
4\,A_2^2+(10\,M_2-M_3)\,A_3^{4/5}L_1^{2/5}
+25\,M_1A_2A_3^{2/5}L_1^{1/5}
)\,L_{1,u}^2-162\,L_1^{-1}L_{1,x}^2
\right.
\\
\fl
&&
\left.
+A_3^{-2}(
8\,A_2^3+135\,A_0A_3^2+6\,(10\,M_2-M_3)\,A_2 A_3^{4/5}L_1^{2/5}
+75\,M_1A_2^2A_3^{2/5}L_1^{1/5}
\right.
\\
\fl
&&
\left.
-M_4A_3^{6/5}L_1^{3/5}
)\,L_{1,u}
-270\,L_1A_{1,x}
+675\,L_1A_{0,u}
+2\,A_2^4A_3^{-2}L_1
+25\,M_1A_2^3A_3^{-8/5}L_1^{6/5}
\right.
\\
\fl
&&
\left.
+3\,(10\,M_2-M_3)\,A_2^2A_3^{-6/5}L_1^{7/5}
+A_2\,L_1(M_4A_3^{-4/5}L_1^{3/5}+540\,A_0)
+162\,A_1^2L_1
\right),
\\
\fl
M_6 &=&-A_3^{-6/5}L_1^{-18/5}\,\left(
2\,L_{1,u}^3
+3\,(2\,A_2L_1+5\,M_1A_3^{2/5}L_1^{6/5})\,L_{1,u}^2
+6\,L_1^2\,(A_2^2+M_1A_2A_3^{2/5}L_1^{1/5}
\right.
\\
\fl
&&
\left.
+M_2A_3^{4/5}L_1^{2/5})\,L_{1,u}
+9\,A_3L_1^3\,(A_{1,u}-2\,A_{2,x})
+2\,L_1^3(A_2^3-27\,A_0A_3^2)
+15\,M_1A_2^2A_3^{2/5}L_1^{16/5}
\right.
\\
\fl
&&
\left.
+6\,M_2A_2A_3^{4/5}L_1^{17/5}
\right).
\end{eqnarray*}
The invariant derivatives
$
\mathbb{D}_1 =
A_3^{1/5}L_1^{-2/5}\,\case{\partial}{\partial x}
-A_3^{-4/5}L_1^{-7/5}(L_{1,u}+A_2L_1)\,\case{\partial}{\partial u}
$
and
$
\mathbb{D}_2 = A_3^{-2/5}L_1^{-1/5}\,\frac{\partial}{\partial u}
$
are defined by the requirement that
$dY = \case{1}{3}\,(\mathbb{D}_1(Y)+P_5\,\mathbb{D}_2(Y))\,\omega_1+\mathbb{D}_2(Y)\,\omega_2$
holds for any function $Y(x,u)$.
The second order classifying manifold associated with forms  $\boldsymbol{\omega}$
in the case $\mathscr{B}_{21}$
has the form
\begin{equation}
\mathcal{C}^{(2)}_{\mathscr{B}_{21}}(\boldsymbol{\omega},\mathcal{V}) =
\left\{
\mathbb{D}_1^i\mathbb{D}_2^j(M_m(x,u))\,\,\vert\,\,
0\le i+j \le 2, \,\, 1\le m\le 6, \,\, (x,u) \in \mathcal{V}
\right\},
\label{classifying_set_B_21}
\end{equation}
where $\mathcal{V} \subset J^0(\pi)$ is an open subset such that   $A_3\not = 0$ in all its points.

\vskip 5 pt

The case $\mathscr{B}_{22}$ is defined by the requirement $A_3=0$. This case is separated on two
subcases: subcase $\mathscr{B}_{221}$ corresponds to the condition $N_0 = L_{1,u}+A_1L_1 \not =0$, while
subcase $\mathscr{B}_{222}$ is defined by $N_0\equiv 0$.

\vskip 5 pt

In the case $\mathscr{B}_{221}$
normalization (\ref{second_normalization}) has the form
$b_1 = L_1\,\left|N_0 \right|^{-1/2}$.
This gives the structure equations
\begin{eqnarray*}
d\omega_1 &=& N_1\,\omega_1\wedge\omega_2,
\\
d\omega_2  &=&
\left(
\case{2}{15}\,(5\,N_2+6)\,P_6+N_2
\right)\,\omega_1 \wedge \omega_2
+\omega_1 \wedge \omega_3,
\\
d\omega_3 &=&
\left(
\case{2}{75}\,(5\,N_1+6)\,P_6^2+\case{1}{5}\,(2\,N_2+3\,N_4)\,P_6+N_3
\right)\,\omega_1 \wedge \omega_2
+\case{3}{5}\,(5\,N_1+4)\,\omega_2\wedge\omega_3
\\
&&
+\case{1}{10}\,(2\,(5\,N_1+6)\,P_6+15\,N_2)\,\omega_1\wedge\omega_3,
\end{eqnarray*}
where for the invariant
$P_6 = N_0^3\,L_1^{-1}\,u_x$
the following equation
\[
dP_6 = \case{3}{5}\,((5\,N_1+4)\,P_6-N_5)\,\omega_2
+3\,\omega_3
+\left(
(N_1+1)\,P_6^2+\case{1}{10}\,(15\,N_2-2\,N_5)\,P_6+N_6
\right)\,\omega_1
\]
holds, while the other invariants
\begin{eqnarray*}
N_1 &=& 3\,N_0^{-3}\,\left(L_1\,N_{0,u} + A_2L_1N_0\right) - 3,
\\
N_2 &=& 2\,L_1^{-1}\,(N_{0,x}-A_1N_0)-\case{6}{5}\,N_5,
\\
N_3 &=& \case{1}{5}\,N_0L_1^{-1}\,N_{5,x}+N_0^2L_1^{-1}A_{0,u}
+\case{1}{5}\,N_0^2L_1^{-3}A_0\,\left(N_0^2-5\,A_2L_1\right)
+\case{1}{25}\,N_5^2-\case{1}{10}\,N_2N_5,
\\
N_4 &=& N_0^{-1}\,\left(A_{1,u}-2\,A_{2,x}\right),
\\
N_5 &=& N_0L_1^{-2}\,\left(L_{1,x}-2A_1L_1\right),
\\
N_6 &=& A_0N_0^4L_1^{-3}
\end{eqnarray*}
depend on $x$ and  $u$, the invariant derivatives
$
\mathbb{D}_1 = N_0L_1^{-1}\,\frac{\partial}{\partial x}
$
and
$
\mathbb{D}_2 = N_0^{-2}L_1\,\frac{\partial}{\partial u}
$
are defined by the equation
$dY = \left(\mathbb{D}_1(Y)+\mathbb{D}_2(Y)\,P_6\right)\,\omega_1 +3\,\mathbb{D}_2(Y)\,\omega_2$.
The second order classifying manifold associated with forms  $\boldsymbol{\omega}$
in the case $\mathscr{B}_{221}$
has the form
\begin{equation}
\mathcal{C}^{(2)}_{\mathscr{B}_{221}}(\boldsymbol{\omega},\mathcal{V}) =
\left\{
\mathbb{D}_1^i\mathbb{D}_2^j(N_m(x,u))\,\,\vert\,\,
0\le i+j \le 2, \,\, 1\le m\le 6, \,\, (x,u) \in \mathcal{V}
\right\},
\label{classifying_set_B_221}
\end{equation}
where $\mathcal{V} \subset J^0(\pi)$ is an open subset such that $A_3\equiv 0$ and $N_0\not = 0$ in all its points.

\vskip 5 pt

In the case  $\mathscr{B}_{222}$ which is defined by the conditions  $A_3 \equiv 0$ and
$N_0 \equiv 0$ we get the structure equations
\begin{eqnarray*}
d\omega_1 &=& \eta_1\wedge \omega_1,
\\
d\omega_2 &=& -2\,\eta_1\wedge\omega_2 +\omega_1\wedge \omega_3,
\\
d\omega_3 &=&
-3\,\eta_1 \wedge \omega_1
+\case{1}{5}\,b_1^{-2}\,\left(3\,(A_{1,u}-2\,A_{2,x})\,u_x
+ 5\,Q_0\,L_1^{-2}\right)\,\omega_1 \wedge \omega_2,
\end{eqnarray*}
where we denote
\[
Q_0=\case{1}{5}\,(L_1L_{1,xx}-2A_{1,x}L_1^2)
+\case{1}{25}\,(6\,A_1^2L_1^2-A_1L_1L_{1,x}-6\,L_{1,x}^2)
+L_1^2\,(A_{0,u}-A_0A_2).
\]

We consider two cases in the further analysis, the case $\mathscr{B}_{2221}$
is defined by the condition
$V_0=A_{1,u}-2\,A_{2,x} \not = 0$, while the case
$\mathscr{B}_{2222}$ is defined by the condition  $A_{1,u}-2\,A_{2,x} \equiv 0$.

\vskip 5 pt

In the case  $\mathscr{B}_{2221}$
the normalization  $b_1 = \left|15\,(A_{1,u}-2\,A_{2,x})\,u_x+25\,Q_0\,L_1^{-2}\right|^{1/2}$
gives the struc\-tu\-re equation
$
d\omega_1 = -\case{225}{2}\,P_7^{1/2}\,\omega_1\wedge\omega_2+\dots\,\omega_1\wedge\omega_3,
$
with the invariant
$
P_7 = \case{1}{5}\,V_0^2L_1^{-4}\,(3\,V_0L_1^2\,u_x+5\,Q_0).
$
We have
\begin{eqnarray*}
\omega_1 &=& \case{1}{5}\,P_7^{1/2}\,L_1V_0^{-1}\,dx,
\\
\omega_2 &=& \case{1}{75}\,P_7^{-1}\,L_1^{-1}V_0\,(du-u_x\,dx)
\\
\omega_3 &=&\case{1}{375}\,P_7^{-3/2}\,
L_1^{-3}V_0^3\,
\left(L_1\,du_x - \case{1}{5}\,(5\,A_2L_1\,u_x +2\,A_1L_1-L_{1,x})\,du
\right.
\\
&&\left.
-\case{1}{5}\,((3\,A_1L_1+L_{1,x})\,u_x+5\,A_0L_1)\,dx
\right).
\end{eqnarray*}
Since forms $\boldsymbol{\omega}$ and function $P_7$ are invariant,
we can divide, without loss of generality, the right-hand sides of forms $\omega_1$, $\omega_2$, and $\omega_3$
by  $\case{1}{5}\,P_7^{1/2}$, $\case{1}{75}\,P_7^{-1}$, and $\case{1}{375}\,P_7^{-3/2}$, respectively.
The obtained forms satisfy the structure equations
\begin{eqnarray*}
d\omega_1 &=& 0
\\
d\omega_2 &=& Q_1\,\omega_1\wedge \omega_2 +\omega_1\wedge \omega_3,
\\
d\omega_3 &=& P_7\,\omega_1\wedge \omega_2 +\case{3}{2}\,Q_1\,\omega_1\wedge \omega_3,
\end{eqnarray*}
while the differential of the invariant  $P_7$ acquires the form
\[
dP_7 = \case{1}{18}\,((31\,Q_1+10)\,P_7+18\,Q_2)\,\omega_1
+\case{1}{15}\,(2\,Q_1+5)\,\omega_2+\case{3}{5}\,\omega_3,
\]
and invariants
\begin{eqnarray*}
\fl
Q_1&=&2\,V_{0,x}L_1^{-1}-\case{2}{5}\,V_0L_1^{-2}\,(3\,L_{1,x}-A_1L_1),
\\
\fl
Q_2 &=& \case{3}{5}\,A_0V_0^4L_1^{-3}
+\case{1}{5}\,V_0^3L_1^{-6}\,(5\,L_1\,Q_{0,x}-2\,Q_0\,(A_1L_1+7\,L_{1,x}))
-\case{1}{18}\,Q_0V_0^2L_1^{-4}\,(13\,Q_1+10)
\end{eqnarray*}
depend on $x$ and $u$. Then equation
$
dY = \left(\mathbb{D}_1(Y)+\case{3}{5}\,P_7\,\mathbb{D}_2(Y)\right)\,\omega_1
+\mathbb{D}_2(Y)\,\omega_2
$
defines the invariant derivatives
$
\mathbb{D}_1 = V_0L_1^{-1}\,\frac{\partial}{\partial x}
-\case{5}{3}\,V_0L_1^{-3}\,\frac{\partial}{\partial u}
$
and
$
\mathbb{D}_2 = V_0^{-3}L_1\,\frac{\partial}{\partial u}.
$
The second order classifying manifold associated with forms  $\boldsymbol{\omega}$
in the case $\mathscr{B}_{2221}$
has the form
\begin{equation}
\mathcal{C}^{(2)}_{\mathscr{B}_{2221}}(\boldsymbol{\omega},\mathcal{V}) =
\left\{
\mathbb{D}_1^i\mathbb{D}_2^j(Q_m(x,u))\,\,\vert\,\,
0\le i+j \le 2, \,\, 1\le m\le 2, \,\, (x,u) \in \mathcal{V}
\right\},
\label{classifying_set_B_2221}
\end{equation}
where $\mathcal{V} \subset J^0(\pi)$ is an open subset such that $A_3\equiv N_0\equiv 0$ and
$V_0\not = 0$  in all its points.

\vskip 5 pt

In the case $\mathscr{B}_{2222}$, which is defined by the requirement $A_{1,u}-2\,A_{2,x} \equiv 0$,
there exists a fucntion  $B(x,u)$ such that $A_1=2\,B_x$, $A_2=B_u$, that is, equation
(\ref{LLT}) in this case has the form
\begin{equation}
u_{xx}= B_u\,u_x^2+2\,B_x\,u_x+A_0.
\label{BabichBordag}
\end{equation}
We use the following result of \cite{BabichBordag1999}:

\vskip 5 pt
\noindent
{\sc Lemma 3}:
{\it
Each equation (\ref{BabichBordag}) can be mapped into equation of the form
\begin{equation}
u_{xx} = A_0(x,u),
\label{generalizedEmdenFowler}
\end{equation}
by means of a transformation from
}
$\mathrm{Cont}_0(J^2(\pi))$.

\vskip 5 pt
\noindent
{\it Proof}:
let  $U=U(x,u)$ be a function such that $U_u = \exp(-B)$.
Then the direct check shows that the change of variables $\widetilde{x}=x$, $\widetilde{u}=U(x,u)$
transforms equation (\ref{BabichBordag}) into equation (\ref{generalizedEmdenFowler}).
\hfill {\sc qed}

\vskip 5 pt

For equation (\ref{generalizedEmdenFowler}) the structure equations acquire the form
\begin{eqnarray*}
d\omega_1 &=& \eta_1\wedge \omega_1,
\\
d\omega_2 &=& -2\,\eta_1\wedge\omega_2 +\omega_1\wedge \omega_3,
\\
d\omega_3 &=&
-3\,\eta_1 \wedge \omega_1
+\dots\,\omega_1 \wedge \omega_2
+\case{4}{5}\,b_1^2\,A_{0,uuu}\,A_{0,uu}^{-2}\,\omega_2\wedge\omega_3.
\end{eqnarray*}
The further analysis dependes on whether condition $A_{0,uuu}\not= 0$ holds.

\vskip 5 pt
In the case  $\mathscr{B}_{22221}$, which is defined by this condition, the normalization
$b_1= A_{0,uu}\,\left| A_{0,uuu} \right|^{-1/2}$ gives the structure equations
\begin{eqnarray*}
d\omega_1 &=& R_1\,\omega_1\wedge\omega_2,
\\
d\omega_2 &=& \left(\case{2}{5}\,(5\,R_1+2)\,P_8+R_2\right)\,\omega_1\wedge\omega_2
+\omega_1\wedge\omega_3,
\\
d\omega_3 &=& \left(\case{2}{25}\,(5\,R_1+2)\,P_8^2+\case{2}{5}\,R_2\,P_8+R_3\right)\,\omega_1\wedge\omega_2
+\case{3}{10}\,(2\,(5\,R_1+2)\,P_8+5\,R_2)\,\omega_1\wedge \omega_3
\\
&&
+\left(3\,R_1+\case{4}{5}\right)\,\omega_2\wedge \omega_3,
\end{eqnarray*}
where the invariant $P_8 = \left| A_{0,uuu} \right|^{3/2}A_{0,uu}^{-2}\,u_x$ depends on $x$, $u$, $u_x$, while
\[
dP_8 = \left((3\,R_1+1)\,P_8^2+\left(R_4+\case{3}{2}\,R_2\right)\,P_8+R_5\right)\,\omega_1
+\left(\left(3\,R_1+\case{4}{5}\right)\,P_8+R_4\right)\,\omega_2
+\omega_3,
\]
and the invariants
\begin{eqnarray*}
R_1&=& \case{1}{2}\,A_{0,uuuu}A_{0,uu}A_{0,uuu}^{-2}-1,
\\
R_2 &=& (5\,R_1+1)\,R_4-5\,R_{4,u}A_{0,uu}A_{0,uuu}^{-1},
\\
R_3 &=& A_{0,u}A_{0,uuu}A_{0,uu}^{-2}-R_{4,x}\left|A_{0,uuu}\right|^{1/2}A_{0,uu}^{-1}
+\case{1}{5}\,R_5+\case{1}{2}\,R_4\,(R_2+2\,R_4),
\\
R_4 &=& -\case{1}{5}\,A_{0,xuu}\left| A_{0,uuu} \right|^{1/2}A_{0,uu}^{-2},
\\
R_5 &=& A_0A_{0,uuu}^2A_{0,uu}^{-3}
\end{eqnarray*}
depend on $x$, $u$. The identity
$
dY = (\mathbb{D}_1(Y)+\mathbb{D}_2(Y)\,P_8)\,\omega_1
+\mathbb{D}_2(Y)\,\omega_2
$
defines the invariant derivatives
$
\mathbb{D}_1 = \left| A_{0,uuu} \right|^{1/2}\,A_{0,uu}^{-2}\,\frac{\partial}{\partial x}
$
and
$
\mathbb{D}_2 = A_{0,uuu}^{-1}A_{0,uu}\,\frac{\partial}{\partial u}.
$
The second order classifying manifold associated with forms  $\boldsymbol{\omega}$
in the case $\mathscr{B}_{22221}$
has the form
\begin{equation}
\mathcal{C}^{(2)}_{\mathscr{B}_{22221}}(\boldsymbol{\omega},\mathcal{V}) =
\left\{
\mathbb{D}_1^i\mathbb{D}_2^j(R_m(x,u))\,\,\vert\,\,
0\le i+j \le 2, \,\, 1\le m\le 5, \,\, (x,u) \in \mathcal{V}
\right\},
\label{classifying_set_B_22221}
\end{equation}
where $\mathcal{V} \subset J^0(\pi)$ is an open subset sucht that  $A_{0,uuu}\not = 0$ in all its points.

\vskip 5 pt

In the case $\mathscr{B}_{22222}$, when   $A_{0,uuu} \equiv 0$,
the form of equation
(\ref{generalizedEmdenFowler}) can be specified.
\vskip 7 pt

\noindent
{\sc Lemma 4}:
{\it
Each non-linearizable equation (\ref{generalizedEmdenFowler}) with $A_{0,uuu} \equiv 0$
can be mapped into equation of the form
}
\begin{equation}
u_{xx} = u^2+a_0(x)
\label{reduced_generalizedEmdenFowler}
\end{equation}
{\it
by means of a transformation from the pseudogroup
(\ref{point_transform}).}

\vskip 5 pt
\noindent
{\it Proof}: if $A_{0,uuu} \equiv 0$, then a non-linearizable equation (\ref{generalizedEmdenFowler})
has the form $u_{xx}= a_2(x)\,u^2+a_1(x)\,u+a_0(x)$ with $a_2\not = 0$.
The change of variables $\widetilde{x}=\varphi(x)$,
$\widetilde{u} = (a_2(x))^{2/5}\,u+b_0(x)$, where $\varphi(x)$ is a function such that
$\varphi_x = a_2^{2/5}$, while
$b_0 = \case{1}{50}\,a_2^{-14/5}\,(5\,a_2a_{2,xx}-6\,a_{2,x}^2+25\,a_1a_2^2)$,
maps this equation into equation
$
\widetilde{u}_{\widetilde{x}\widetilde{x}} = \widetilde{u}^2+\widetilde{a}_0(\widetilde{x})
$
with
$
\widetilde{a}_0(\varphi(x)) = (a_2(x))^{-4/5}(b_{0,xx}(x)+a_0(x))
-\case{2}{5}\,(a_2(x))^{-9/5}a_{2,x}(x)b_{0,x}(x)-(b_0(x))^2
$.
\hfill {\sc qed}

\vskip 5 pt

For equation (\ref{reduced_generalizedEmdenFowler}) the structure equations are of the form
\begin{eqnarray*}
d\omega_1 &=& -\case{1}{2}\,\omega_1\wedge\omega_2,
\\
d\omega_2 &=& -P_9\,\omega_1\wedge \omega_2 +\omega_1\wedge\omega_3,
\\
d\omega_3 &=& 2\,\omega_1\wedge\omega_2-\case{3}{2}\,P_9\,\omega_1\wedge\omega_3
-\case{3}{2}\,\omega_2\wedge\omega_3,
\end{eqnarray*}
where $P_9 = u_x\,u^{-3/2}$.

\vskip 5 pt
In the case $\mathscr{B}_{222221}$ such that $a_0\equiv 0$,
that is, in the case when equation
(\ref{reduced_generalizedEmdenFowler}) is $u_{xx}=u^2$, we have
\[
dP_9 = \left(1-\case{3}{2}\,P_9^2\right)\,\omega_1 -\case{3}{2}\,P_9\,\omega_2.
\]

In the case $\mathscr{B}_{222222}$ such that $a_0\not = 0$ we obtain
\[
dP_9 = \left(S_1+1-\case{3}{2}\,P_9^2\right)\,\omega_1 -\case{3}{2}\,P_9\,\omega_2,
\]
where $S_1 = a_0\,u^{-2}$. Consider the subcase $\mathscr{B}_{2222221}$ such that
$a_{0,x} \equiv 0$. In this case we get
\[
d S_1= -2\,S_1\,(P_9\,\omega_1+\omega_2).
\]
Therefore all the equations (\ref{reduced_generalizedEmdenFowler})
in this case have the same structure equations and are equivalent to each other, in particular,
all of them are equivalent to  equation $u_{xx}= u^2+1$.

\vskip 5 pt
Finally, in the last subcase  $\mathscr{B}_{2222222}$ such that $a_{0,x} \not =0$ we obtain
\[
dS_1 = S_1\,(T_1S_1^{1/4}-2\,P_9)\,\omega_1-2\,S_1\,\omega_2,
\]
where  $T_1 = a_{0,x}a_0^{-5/4}$,
and
\[
dT_1 = \left(T_2 -\case{5}{4}\,S_2^2\right)\,\omega_1,
\]
where $T_2 = a_{0,xx}a_0^{-3/2}$.
Therefore
the first order classifying manifold associated with forms  $\boldsymbol{\omega}$
in the case $\mathscr{B}_{2222222}$ for equation (\ref{reduced_generalizedEmdenFowler})
can be taken in the form
\begin{equation}
\mathcal{C}^{(1)}_{\mathscr{B}_{2222222}}(\boldsymbol{\omega},\mathcal{I}) =
\left\{
(a_{0,x}(x)\,(a_0(x))^{-5/4},a_{0,xx}(x)(a_0(x))^{-3/2})\,\,\vert\,\,
x \in \mathcal{I}
\right\},
\label{classifying_set_B_2222222}
\end{equation}
where $\mathcal{I} \subset \mathbb{R}$ is an open interval such that $a_{0,x}\not = 0$ in all its points.

\vskip 7 pt

Combining the results of the above computations and applying Theorem 14.24 from \cite{Olver1995},
we get the following theorem.

\vskip 7 pt
\noindent
{\sc Theorem 2}:
{\it
Each equation (\ref{LLT})
can be transformed to an equation from one of the invariant subclasses
$\mathscr{B}_{1}$, ... ,  $\mathscr{B}_{2222222}$ by a diffeomorphsim for the pseudogroup of
point transformations (\ref{point_transform}).

Subclass  $\mathscr{B}_{1}$ contains equations  (\ref{LLT}) with $L_1\equiv L_2\equiv 0$.
These equations are locally equivalent to equation $u_{xx}=0$.

Equations (\ref{LLT}) such that one of the functions $L_1$ or $L_2$ is not equal to zero,
can be mapped to an equation with $L_1\not=0$, $L_2\equiv 0$.
These equations are divided on invariant subclasses
$\mathscr{B}_{21}$, ... , $\mathscr{B}_{2222222}$.

Subclass  $\mathscr{B}_{21}$  contains equations (\ref{LLT}) with $L_1\not =0$, $L_2\equiv 0$, $A_3\not =0$.
Two equations from this subclass are locally equivalent with respect to the pseudogroup (\ref{point_transform})
whenever their classifying manifolds
(\ref{classifying_set_B_21}) locally overlap.

Subclass  $\mathscr{B}_{221}$ contains equations (\ref{LLT}) with
$L_1\not =0$, $L_2\equiv A_3\equiv0$,  $N_0\not = 0$.
Two equations from this subclass are locally equivalent with respect to the pseudogroup (\ref{point_transform})
whenever their classifying manifolds
(\ref{classifying_set_B_221}) locally overlap.

Subclass  $\mathscr{B}_{2221}$  contains equations (\ref{LLT}) with
$L_1\not =0$, $L_2\equiv A_3\equiv N_0 \equiv 0$,  $V_0\not = 0$.
Two equations from this subclass are locally equivalent with respect to the pseudogroup (\ref{point_transform})
whenever their classifying manifolds
(\ref{classifying_set_B_2221}) locally overlap.

Equations (\ref{LLT}) with $L_1\not =0$, $L_2\equiv A_3\equiv N_0 \equiv V_0 \equiv 0$
can be mapped into equations of the form   (\ref{generalizedEmdenFowler}).

Subclass  $\mathscr{B}_{22221}$ contains equations (\ref{generalizedEmdenFowler}) with
$A_{0,uuu}\not = 0$.
Two equations from this subclass are locally equivalent with respect to the pseudogroup (\ref{point_transform})
whenever their classifying manifolds
(\ref{classifying_set_B_22221}) locally overlap.

Equations (\ref{generalizedEmdenFowler}) with $A_{0,uuu}\equiv 0$
can be mapped into equations of the form  (\ref{reduced_generalizedEmdenFowler}).

Subclass  $\mathscr{B}_{222221}$  contains one equation $u_{xx}=u^2$.

Subclass  $\mathscr{B}_{2222221}$  contains equations$u_{xx}=u^2+\alpha$, $\alpha = \mathrm{const}$,
$\alpha\not=0$. All these equations are equivalent to each other, in particular, all of them are equivalent to
equation $u_{xx}=u^2+1$.

Subclass $\mathscr{B}_{2222222}$ contains equations (\ref{reduced_generalizedEmdenFowler}) with
$a_{0,x}\not = 0$.
Two equations from this subclass are locally equivalent with respect to the pseudogroup (\ref{point_transform})
whenever their classifying manifolds
(\ref{classifying_set_B_2222222}) locally overlap.

}

\end{document}